\documentclass[a4paper,12pt]{article}

\usepackage{amsmath,amssymb,amsthm}
\usepackage{bbm, mathtools}
\usepackage{color}
\usepackage[mathletters]{ucs}
\usepackage[utf8x]{inputenc}
\usepackage{hyperref}
\newtheorem{theorem}{Theorem}

\newcommand{\thmref}[1]{Theorem~\ref{#1}}

\newcommand{\xd}{||x||_2}

\def\={\stackrel {\rm def}  {=}}
\def\di={\overset{\text{${\mathcal{D} }$}} =}

\title{A non-uniform Littlewood-Offord inequality}
\author{D. Dzindzalieta$^{1}$, T. Ju\v skevi\v cius$^{2}$}

\begin{document}
\maketitle
\begin{abstract}
Consider a sum $S_n=v_i\varepsilon_1+\cdots+v_n\varepsilon_{n}$, where $(v_i)^{n}_{i=1}$ are non-zero vectors in $\mathbb{R}^{d}$ and $(\varepsilon_i)^{n}_{i=1}$ are independent Rademacher random variables (i.e., $~{\mathbb{P}(\varepsilon_{i}=\pm 1)=~\frac{1}{2}}$). The classical Littlewood-Offord problem asks for the best possible upper bound for $~{\sup_{x}\mathbb{P}(S_n = x)}$.
In this paper we consider a non-uniform version of this problem. Namely, we obtain the optimal bound for $\mathbb{P}(S_n = x)$ in terms of the length of the vector $x\in \mathbb{R}^d$.

\end{abstract}
\section{Introduction}
Let $(\varepsilon_i)^{n}_{i=1}$ be a collection of independent Rademacher random variables and denote their sum by $R_n$. That is, we have $\mathbb{P}(\varepsilon_i=\pm 1)=\frac{1}{2}$. We shall throughout the paper denote by $S_n$ the weighted sum of Rademacher random variables $v_{1}\varepsilon_{1}+\cdots+v_n\varepsilon_{n}$, where weights $v_i\neq 0$ are vectors in $\mathbb{R}^d$ such that $||v_i||_{2}\leq 1$. On each occasion we shall specify which dimension we are working in. Define the quantity $\delta_{n,k}$ to be equal to $0$ if $n+k\in 2\mathbb{Z}$ and $1$ otherwise.

The classical Littlewood-Offord problem asks for the best possible bound for $\mathbb{P}(S_n =~ x)$. It turns out that for all $d$ we have
\begin{theorem}
\begin{equation}\label{erd}
\mathbb{P}(S_n = x)\leq \mathbb{P}(R_n = \pm \delta_{n,0})=\binom{n}{\lfloor\frac{n}{2}\rfloor}/2^n.
\end{equation}
\end{theorem}
The latter bound is clearly optimal and was established by Erd\H{o}s \cite{ELO} in the case $d=1$. Answering a question of Erd\H{o}s, Kleitman \cite{kleitman} extended his result to all $d\geq 1$. \\

Since the introduction of the Littlewood-Offord problem, many variations on the problem have been considered. Answering a question of Erd\H{o}s and Moser \cite{extremal}, Sarkozy and Szemeredi \cite{SS} proved that $\mathbb{P}(S_n = x)=O(n^{-\frac{3}{2}})$ in $d=1$ under the assumption that all $v_i$'s are distinct. The latter research culminated in Stanley's famous proof \cite{Stanley} of the exact bound conjectured by Erd\H{o}s and Moser \cite{extremal}, namely that the choice $v_i=i$ is optimal. Fairly recently, Tao and Vu \cite{TV, TV2} and Nguyen and Vu \cite{VN} investigated inverse Littlewood-Offord problems that are now an important tool in Random matrix theory. Their results can informally be described by saying that if $\mathbb{P}(S_n = x)$ is "large" then $v_i$'s can be covered by a "small" number of generalized arithmetic progressions. Bandeira, Ferber and Kwan \cite{BFK} considered the resilience version of the Littlewood-Offord problem and have formulated some very interesting open questions. Tiep and Vu \cite{Tiep} were the first ones to consider this problem in a non-Abelian setting, namely, for certain matrix groups. Their work has been very recently extended and optimal bounds obtained by Ju\v skevi\v cius and \v Semetulskis in arbitrary groups \cite{JS}.\\

In this paper we shall establish a non-uniform bound for $\mathbb{P}(S_n = x)$ in terms of the length of $x\in \mathbb{R}^d$ and $n$. The main result of the paper is the following.

\begin{theorem}\label{nonuniform}
For arbitrary $d$ and all non-zero $x\in \mathbb{R}^d$ we have
\begin{equation}
\mathbb{P}(S_n =x)\leq \mathbb{P}(R_n = k+\delta_{n,k})=\binom{n}{\lceil \frac{n+k}{2}\rceil}/2^n,
\end{equation}
where $k$ is the upper integer part of $\xd$.
\end{theorem}

It is not difficult to see that the latter bound is optimal. Equality is achieved by the sum $S_n=(\frac{||x||_2}{k+\delta_{n,k}}R_n,0,\ldots,0)$. \thmref{nonuniform} in $d=1$ was proved in \cite[Theorem 1.1]{DJS}. In his proof of $\eqref{erd}$ Erd\H{o}s used a result from extremal set theory - Sperner's Theorem. The inequality in \thmref{nonuniform} for $d=1$ is proved in a similar fashion, but another tool is needed - Milner's Theorem on the size of intersecting antichains. We have included the proof of the case $d=1$ for the sake of completeness.  It turns out that the one-dimensional result can be then extended to all dimensions by induction together with Kleitman's bound \eqref{erd}.

Kleitman's result tells us that for all $x\in \mathbb{R}^d$ we have $\mathbb{P}(S_n=x)=O(n^{-\frac{1}{2}})$. ~{\thmref{nonuniform}} gives us more detailed information - it tells us that the probability in question is exponentially small for $||x||_2$ much larger than $\sqrt{n}$. To give a crude bound one can use Hoeffding's inequality that gives us
$\mathbb{P}(S_n=x)\leq \exp (-||x||_{2}^2/(2n)$). Of course, one can use more detailed analysis using Stirling's approximation to obtain sharp asymptotic expressions.\\
\textbf{Remark}. \thmref{nonuniform} easily extends to arbitrary symmetric random variables $X_i$ such that $||X_i||_2 \leq 1$ and $\mathbb{P}(X_i=0)=0$ by conditioning on the norm of the variables $X_i$ and using the statement of the theorem for two-point distributions.\\

The requirement that $v_i \neq 0$ is indeed essential and we shall illustrate it with the following result that follows easily from the one-dimensional case and that actually was our first result in trying to generalize the results from \cite{DJS} to high dimensions.

\begin{theorem}\label{zeroweights}

Let $X_1,\ldots, X_n$ be independent symmetric random variables in $\mathbb{R^d}$ such that  $||X_i||_{2}\leq 1$. Then for all non-zero $x\in \mathbb{R}^d$ we have

\begin{equation}
\sup_{n\geq 1}\mathbb{P}\left\{X_{1}+\cdots+X_{n} = x\right\}= \mathbb{P}\left\{R_{k^2} = k\right\},
\end{equation}
where $k$ is the upper integer part of $||x||_{2}$.
\end{theorem}

The latter result shows that even in the case of $X_i=v_i\varepsilon_i$ the probability $\mathbb{P}(S_n=x)$ can be bounded away from zero as $n\rightarrow \infty$ if we allow zero weights. In contrast to this situation, \thmref{nonuniform} tells us that $\mathbb{P}(S_n=x)=O(n^{-\frac{1}{2}})$ if $v_i\neq 0$.\\

Finally, we address the situation not covered by \thmref{nonuniform}, namely, the case $x=0$. Note that for $n\in 2\mathbb{Z}$ Kleitman's result gives the best possible bound for $\mathbb{P}(S_n=0)$. The situation for $n\in 2\mathbb{Z}+1$ is subtly different due to parity issues.

\begin{theorem}\label{zerocase}
Under the assumptions of \thmref{nonuniform} and $n\in 2\mathbb{Z}+1$ we have
$$\mathbb{P}(S_n=0)\leq \mathbb{P}\left(\frac{1}{2}R_{n-1}+\varepsilon_{n}=0\right).$$

\end{theorem}

\section{Proofs}

We shall first provide the proof of \thmref{zeroweights} for the reason that it is is very easy to deduce it from the case $d=1$ proved in \cite[Theorem 1.2]{DJS}.

\textit{Proof of \thmref{zeroweights}.}\, Due to the symmetry of  $X_i$'s the sums $~{X_1+\cdots+X_n}$ and $X_1\varepsilon_1+\cdots+X_n\varepsilon_n$ have the same distribution. Condition on the sequence $X_1,\ldots,X_n$. Then $X_1\varepsilon_1+\cdots+X_n\varepsilon_n$ is distributed as $v_1\varepsilon_1+\cdots+v_n\varepsilon_n$, where $v_i=X_i$ and $\varepsilon_i$'s are independent of the $X_i$'s. We have that $||X_i||_2=||v_i||_2\leq 1$. First note that we can assume that $x$ is one-dimensional, that is, we can assume that $x=(\xd,0,\ldots,0)$. This can be achieved by changing the basis of $\mathbb{R}^d$ by an orthogonal transformation that also preserved the lengths of the vectors $v_i$. Write $v^{(1)}$ for the first coordinate of the vector $v$. We have
$$\mathbb{P}(S_n = x)\leq \mathbb{P}(v_1^{(1)}\varepsilon_1+\cdots+v_n^{(1)}\varepsilon_n = \xd)$$
and the desired result then follows from Theorem 1.2 from \cite{DJS}.

\textit{Proof of \thmref{nonuniform} in the case $d=1$.}\, Let $v_1,\ldots,v_n$ be non-zero real numbers satisfying $|v_i|\leq 1$. The distribution of $S_n=v_{1}\varepsilon_{1}+\cdots+v_{n}\varepsilon_{n}$ is unchanged if we replace $v_i$ by $-v_i$. Therefore we can assume that $v_i>0$ and $x>0$. To obtain the desired inequality we shall use a result in extremal combinatorics due to Milner \cite{milner}. We shall say that a family of subsets $\mathcal{F}$ of $[n]$ is an \textit{antichain} if for all $A,B\in \mathcal{F}$ we have $A\not\subset B$ and $k$-\textit{intersecting} if for all $A,B\in \mathcal{F}$ we have $|A\cap B|\geq k$. Milner \cite{milner} proved that if a family of subsets $\mathcal{F}$ of $[n]$ is a $k$-intersecting antichain, then
$$ |{\cal F}| \leq \binom{n}{t}, \qquad t = \left\lceil \frac {n+k} 2 \right\rceil.$$
For a sum $S_n=v_{1}\varepsilon_{1}+\cdots+v_{n}\varepsilon_{n}$ define the family of subsets $\mathcal{F}_x =\{ A\subset [n]:\, \sum_{i\in A}v_i-\sum_{A^{c}}v_i=x\}$. For notational convenience let us denote by $\sigma_A$ the sum $\sum_{i\in A}v_i$ and $s_A=\sigma_A-\sigma_{A^c}$. It is easy to see that $\mathcal{F}_x$ is an antichain. Indeed, for any distinct subsets $A, B$ of $[n]$ we that $A\subset B$ implies that $s_A<s_B$ as $v_i>0$ and we are done. Let us now show that $\mathcal{F}_x$ is $k$-intersecting with $k=\lceil x\rceil$. Assume that $A,B \in \mathcal{F}_x$ and that $|A\cap B|<k$. We then have
  \begin{equation}\label{one}
    s_A = \sigma_A - \sigma_{A^{c}} =  (\sigma_{A \cap B} - \sigma_{A^c \cap B^{c}}) + (\sigma_{A\cap B^{c}} - \sigma_{A^{c}\cap B})                                                                                                                        \end{equation}
  and
  \begin{equation}\label{two}
    s_B = \sigma_B - \sigma_{B^{c}} = (\sigma_{A \cap B} - \sigma_{A^{c}\cap B^{c}}) - (\sigma_{A \cap B^{c}} - \sigma_{A^{c}\cap B}).
  \end{equation}
  Since
  $$\sigma_{A \cap B} - \sigma_{A^{c}\cap B^{c}} \leq \sigma_{A \cap B} \leq |A \cap B| \leq k -1 < x,$$
  from \eqref{one} and \eqref{two} we get
  $$\min\{s_A, s_B\} < x,$$
  which contradicts the fact $s_A= s_B = x$.\\
To complete the proof we just note that
$$\mathbb{P}(S_n=x)=|\mathcal{F}_x|/2^{n}\leq \binom{n}{\lceil \frac{n+k}{2}\rceil}/2^n=\mathbb{P}(R_n = k+\delta_{n,k}).$$

\textit{Proof of \thmref{nonuniform}.}\,For a vector $v\in \mathbb{R}^d$ we denote its $j$-th coordinate by $v^{(j)}$.  Without loss of generality we can assume that $~{x=(\xd,0,\ldots, 0)}$ since if that is not the case, we can change the coordinate system as in the proof of \thmref{zeroweights} and achieve this while we keep the lengths of the vectors $v_i$ unchanged .
 Let $m$ be the number $v_i$'s with non-zero first coordinate $v_{i}^{(1)}$. Without loss of generality we assume that these are the first $m$ vectors $v_i$. If $m=0$ the problem reduces to $(d-1)$-dimensions. We shall from now on assume that $m\geq 1$.

Let us write $\mathcal{E}$ for the collection of random variables $\varepsilon_i$ with indices $i\leq m$.  Given a realization of $\mathcal{E}$ write $s_{\mathcal{E}}=\sum_{i=0}^m v_i^{\ast}\varepsilon_i$, where $v_i^{\ast}=(0,v_i^{(2)},\ldots,v_i^{(d)})$. Denote by $k$ an upper integer part of $\xd$. We have 
\begin{eqnarray}
&&\mathbb{P}(v_1\varepsilon_1+\cdots+ v_n\varepsilon_n = x)=\mathbb{E}_{\mathcal{E}}\mathbb{P}(\sum_{i=0}^m v_i\varepsilon_i+\sum_{i=m+1}^n v_i\varepsilon_i=x|\mathcal{E}) \label{pirma}\\
&&=\mathbb{E}_{\mathcal{E}}\mathbb{P}(\sum_{i=1}^m v^{(1)}_i\varepsilon_i=\xd,\sum_{i=m+1}^nv_i\varepsilon_i=-s_{\mathcal{E}}|\mathcal{E})\label{antra} \\
&&= \mathbb{P}(\sum_{i=1}^mv_i^{(1)}\varepsilon_i=\xd )\mathbb{E}_{\mathcal{E}}\mathbb{P}(\sum_{i=m+1}^nv_i\varepsilon_i=-s_{\mathcal{E}}|\mathcal{E}) \label{trecia}\\
&&\leq  \mathbb{P}(\sum_{i=1}^{m}\varepsilon_i=k+\delta_{m, k})\mathbb{E}_{\mathcal{E}} \mathbb{P}(\sum^{n}_{i=m+1}\varepsilon_i=-\delta_{n-m, 0} )\label{ketvirta}\\
&&=  \mathbb{P}(\sum_{i=1}^{m}\varepsilon_i=k+\delta_{m, k})\mathbb{P}(\sum^{n}_{i=m+1}\varepsilon_i=-\delta_{n-m, 0} )\notag \\
&&=\mathbb{P}(\sum_{i=1}^{m}\varepsilon_{i}=k+\delta_{m, k})\mathbb{P}(\sum^{n}_{i=m+1}\varepsilon_i=(-1)^{m+k}\delta_{n-m, 0} )\notag\\
&&\leq \mathbb{P}(\sum_{i=1}^{m}\varepsilon_{i}+\sum^{n}_{i=m+1}\varepsilon_i=k+\delta_{m, k}+(-1)^{m+k}\delta_{n-m, 0} )\label{septyni}\\
&&=\mathbb{P}(R_n = k+\delta_{n,k}).\label{astuoni}
\end{eqnarray}
The equality \eqref{pirma} follows from the law of total probability. Equality from \eqref{antra} to \eqref{trecia} follows by independence of the events $\{\sum_{i=1}^{m}v^{(1)}_i\varepsilon_i=\xd\}$ and $\{\sum^{n}_{i=m+1}v_i\varepsilon_i=-s_{\mathcal{E}}\}$
conditioned on $\mathcal{E}$. The inequality from \eqref{trecia} to \eqref{ketvirta} follows from the case $d=1$ and \thmref{erd}. Equality $\delta_{m, k}+(-1)^{m+k}\delta_{n-m, 0}=\delta_{n,k}$ follows from the definition of the function $\delta_{n,k}$ justifying the equality from \eqref{septyni} to \eqref{astuoni}.

\textbf{Remark}. The proof of  \thmref{nonuniform} is more involved than the trivial argument giving us \thmref{zeroweights} using the results from \cite{DJS}. The same argument cannot be used here as after the appropriate rotation in the proof we cannot guarantee that the corresponding new weights are non-zero, which is essential.\\

\textit{Proof of \thmref{zerocase}.}\, We can without loss of generality assume that $||v_n||_{2}=1$. Applying \thmref{nonuniform} we obtain
\begin{eqnarray*}
\mathbb{P}(S_n = 0)&=&\frac{1}{2}\mathbb{P}(S_{n-1} = v_n)+\frac{1}{2}\mathbb{P}(S_{n-1} = -v_n)\\
&\leq&\frac{1}{2}\mathbb{P}(R_{n-1} = 2)+\frac{1}{2}\mathbb{P}(R_{n-1} = 2)\\
&=&\frac{1}{2}\mathbb{P}(R_{n-1} = 2)+\frac{1}{2}\mathbb{P}(R_{n-1} = -2)\\
&=&\mathbb{P}(R_{n-1} +2\varepsilon_{n}= 0)=\mathbb{P}(\frac{1}{2}R_{n-1} +\varepsilon_{n}= 0).\\
\end{eqnarray*}
	
\section{Open questions and conjectures}

In their landmark paper Tao and Vu \cite{TV} proved series of inverse Littlewood-Offord type results that are nowadays a crucial tool in studying discrete random matrices. These inverse results in our setting can be vaguely expressed by saying that if $\sup_{v_i,x}\mathbb{P}(S_n=x)$ is large, then the multiset $\{v_i: i\in [n]\}$ has strong additive structure, meaning that most of the $v_i$'s can be covered by a small number of generalized arithmetic progressions. We thus naturally ask:

\textbf{Question 1.} Suppose that for some $x$ the probability $\mathbb{P}(S_n=x)$ is large, is there a corresponding inverse principle?\\

The latter question is vague as stated, but by it we just mean whether analogous results as in \cite{TV} or even stronger results from the subsequent papers \cite{TV2} and \cite{VN} in this case can be established.\\

We strongly suspect that something very similar to Rademacher random variables should also be true for other types of distributions on the integers. Therefore we formulate the following conjecture.

\textbf{Conjecture 1.} Let $U_1\ldots, U_n$ be independent uniform random variables on the arithmetic progression $\{-m+1,-m+3,\ldots,m-3,m-1\}$ with $m\geq 3$. Then for all non-zero $x\in \mathbb{R}^d$ and non-zero $v_i\in \mathbb{R}^d$ with $||v_i||_{2}\leq 1$ we have that for $m\in 2\mathbb{Z}+1$
$$\mathbb{P}(v_{1}U_{1}+\cdots + v_{n}U_{n}=x)\leq \mathbb{P}(U_{1}+\cdots + U_{n}=k)$$
and for $m\in 2\mathbb{Z}$
$$\mathbb{P}(v_{1}U_{1}+\cdots + v_{n}U_{n}=x)\leq \mathbb{P}(U_{1}+\cdots + U_{n}=k+\delta_{n,k})$$
where $k$ is the lower integer part of $||x||_2$.\\
Note that Conjecture $1$ reduces to \thmref{nonuniform} for $m=2$.

The proof of \thmref{nonuniform} relies on the fact that rotations preserve the Euclidean norm. We nevertheless believe that this is just the limitation of our approach and thus conjecture the following.

\textbf{Conjecture 2.} \thmref{nonuniform} remains true if we replace $||\cdot||_2$ by any other norm on $\mathbb{R}^n$.

\bibliographystyle{unsrt}
\bibliography{mybibliography}

\end{document}